\documentclass[11pt]{article}

\usepackage[margin=1in]{geometry}

\usepackage{amsmath}
\usepackage{multicol}

\usepackage{color}

\begin{document}
%
\title{An MISOCP-Based Solution Approach to \\ the Reactive Optimal Power Flow Problem}
%
%
%

\author{Sezen Ece Kayac{\i}k, Burak Kocuk\thanks{
Industrial Engineering Program, Sabanc{\i} University, Istanbul, Turkey (e-mails: ekayacik, burak.kocuk@sabanciuniv.edu).
}
}

\newcommand{\cI}{\mathcal{I}}
\newcommand{\cS}{\mathcal{S}}
\newcommand{\cT}{\mathcal{T}}
\newcommand{\cB}{\mathcal{B}}
\newcommand{\cL}{\mathcal{L}}
\newcommand{\cG}{\mathcal{G}}

\maketitle

\begin{abstract}
In this letter, we present an alternative mixed-integer non-liner programming formulation of the reactive optimal power flow (ROPF) problem. We utilize a mixed-integer second-order cone programming (MISOCP) based approach to find global optimal solutions of the proposed ROPF problem formulation. We strengthen the MISOCP relaxation via the addition of convex envelopes and cutting planes. Computational experiments on challenging test cases show that the  MISOCP-based approach yields promising results compared to a semidefinite programming based approach from the literature.
\end{abstract}


%


\section{Introduction}
%
%
%
%

{T}{he} reactive optimal power flow (ROPF) problem is a variant of the well-known optimal power flow (OPF) problem in which additional discrete decisions, such as shunt susceptance and tap ratio, are considered. Due to the presence of these discrete variables in the ROPF  problem, it can be formulated as a mixed-integer non-linear programming (MINLP) problem. This letter utilizes the recent developments in the OPF problem to propose an efficient way of solving the ROPF problem.


OPF  is one of the most studied problems in the area of power systems  and
 a variety of solution approaches have been proposed in the literature. Local methods such as the interior point method try to solve the OPF problem but they do not provide any assurances of global  optimality. In recent years, convex relaxations of the OPF problem have drawn considerable research interest since the convexity property promises a globally optimal solution under certain conditions. Several approaches have been developed based on convex quadratic, semidefinite programming (SDP), second order cone programming (SOCP) and convex-distflow formulations. The ROPF problem has a similar structure with the OPF problem, except the inclusion of shunt susceptance and tap ratio variables, which are typically modelled as discrete variables. The resulting MINLP problem is difficult so solve and the literature has primarily focused on {various heuristic methods~\cite{capitanescu2010sensitivity}.} The systematic treatment of the ROPF problem is limited to an SDP-based relaxation called \textit{tight-and-cheap relaxation} (TCR) proposed in \cite{bingane2019tight}. 

This letter proposes a  new MINLP  formulation for the ROPF problem along with its mixed-integer second-order cone programming (MISOCP) relaxation and an improved MISOCP relaxation with convex envelopes  and cutting planes. We also test the accuracy and efficiency of our approach with the TCR method from the literature on difficult test cases and  obtain promising results.

\section{Mathematical Model}
\subsection{MINLP Formulation}
\label{sec:MINLP}
Consider a power network $\mathcal{N} = (\mathcal{B}, \mathcal{L})$, where $\mathcal{B}$ and $\mathcal{L}$ denote the set of buses and the set of transmission lines respectively. Let  $\mathcal{G} \subseteq \mathcal{B}$, $\mathcal{S} \subseteq \mathcal{B}$ and $\mathcal{T} \subseteq \mathcal{L}$ respectively denote the set of generators connected to the grid, the buses with a variable shunt susceptance and the lines with a variable tap ratio.
Rest of the parameters are given as follows:
\begin{itemize}
\item
For each bus $i \in \mathcal{B}$;
$p_i^d$ and $q_i^d$ are the real and reactive power load, $\underline V_i$ and $\overline V_i$ are the bounds on the voltage magnitude,  $\delta(i)$ is the set of neighbors and $\{b_{ii}^k : k \in \mathcal{S}_i\}$ is the set of allowable shunt susceptances.
\item
For each generator located at bus $i \in \mathcal{G}$; active and reactive outputs must be in the intervals $ [p_i^{\text{min}}, p_i^{\text{max}}]$ and $[q_i^{\text{min}} , q_i^{\text{max}}]$, and we have $ p_i^{\text{min}} =  p_i^{\text{max}}=q_i^{\text{min}} =  q_i^{\text{max}}=0$ for $i \in \mathcal{B} \setminus \mathcal{G}$.

\item
For each line $(i,j) \in \mathcal{L}$; $G_{ij}$ and $B_{ij}$ are conductance and succeptance, $\{\tau_{ij}^l : l \in \mathcal{T}_{ij}\}$ is the set of allowable tap ratios, $\overline S_{ij}$ is the apparent power flow limit and $\overline\theta_{ij}$ is the bound on the phase angle.
\end{itemize}
\ \  We define the following decision variables:
\begin{itemize}
\item 
For each bus $i \in \mathcal{B}$, $|V_i|$ and $\theta_i$ are the voltage magnitude and phase angle, $b_{ii}$ is the shunt susceptance, $\alpha_{i}^k $ is one if $b_{ii}=b_{ii}^k$ and zero otherwise.

\item
For each generator located at bus $i \in \mathcal{G}$, $p_i^g$ and $q_i^g$ are the real and reactive power output.

\item
For each line $(i,j) \in \mathcal{L}$, $p_{ij}$ and $q_{ij}$ are the real and reactive power flow, $\tau_{ij}$ is the tap ratio, $\beta_{ij}^l$ is one if $\tau_{ij}=\tau_{ij}^l$ and zero otherwise.
\end{itemize}

Then, the ROPF problem can be modeled as the following MINLP:
\noindent
\begin{flalign}\label{eq:obj R_OPF}
\min  \sum_{i\in\mathcal{G}} f(p_i^g) \hspace{6cm} 
\end{flalign}

\begin{equation}\label{eq:real balance}
p_i^g - p_i^d =  \ g_{ii} |V_i|^2 +  \sum_{j\in\delta(i)} p_{ij} \hspace{3.2cm} i \in \mathcal{B}
\end{equation}
\begin{equation} \label{eq:reactive balance}
q_i^g - q_i^d = -b_{ii} |V_i|^2 + \sum_{j\in\delta(i)} q_{ij} \hspace{3cm} i \in \mathcal{B}
\end{equation}
\begin{equation}\label{eq:real flow}\begin{split}
 p_{ij} = G_{ij}(|V_i|/\tau_{ij})^2  + (|V_i|/\tau_{ij}) |V_j| [G_{ij} \cos(\theta_i- \theta_j) \ \   \\  -B_{ij} \sin(\theta_i - \theta_j) ] \hspace{3.1cm} (i,j) \in \mathcal{L}
\end{split}\end{equation}
\begin{equation}\label{eq:reactive flow}\begin{split}
q_{ij} = - B_{ij}(|V_i|/\tau_{ij})^2   - (|V_i|/\tau_{ij}) |V_j| [B_{ij}\cos(\theta_i - \theta_j) \  \\ 
+G_{ij} \sin(\theta_i - \theta_j) ]  \hspace{3.1cm} (i,j) \in \mathcal{L}
\end{split}\end{equation}
\begin{equation} \label{eq:voltage}
   \underline V_i \le |V_i| \le \overline  V_i \hspace{5cm} i \in \mathcal{B}
\end{equation}
\begin{equation} \label{eq:shunt and tap selection}
   \sum_{k\in \mathcal{S}_i}  b_{ii}^k \alpha_{i}^k = b_{ii} \hspace{0.5cm} i \in \mathcal{B}, \ \sum_{l\in \mathcal{T}_{ij}}  \frac{ \beta_{ij}^l}{\tau_{ij}^l}  = \frac{1}{\tau_{ij}} \hspace{0.5cm} (i,j) \in \mathcal{L}
\end{equation}
\noindent
\begin{equation} \label{eq:binary alpha and beta}
   \sum_{k\in \mathcal{S}_i} \alpha_{i}^k = 1 \hspace{1cm} i \in \mathcal{B}, \ \sum_{l\in \mathcal{T}_{ij}} \beta_{ij}^l = 1 \hspace{1cm} (i,j) \in \mathcal{L}
\end{equation}


\noindent
\begin{equation} \label{eq:binary alpha and beta1}
   \alpha_{i}^k \in \{0,1\} \hspace{1.1cm} i \in \mathcal{B}, \ \beta_{ij}^l \in \{0,1\} \hspace{1.1cm} (i,j) \in \mathcal{L}
\end{equation}
\begin{equation}\label{eq:trivial}
 b_{ii}=0 \hspace{1.75cm}  i\not\in\cS,   \hspace{0.1cm} 
 \tau_{ij}=1 \hspace{1.6cm}  (i,j)\not\in\cT   
\end{equation}
\begin{equation} \label{eq:real and reactive gen}
    q_i^{\text{min}}  \le q_i^g \le q_i^{\text{max}} ,  \hspace{1.2cm} p_i^{\text{min}}  \le p_i^g \le p_i^{\text{max}} \hspace{0.8cm} i \in \mathcal{G}
\end{equation}
\begin{equation} \label{eq:flow and angle limit}
    p_{ij}^2+q_{ij}^2  \le \overline S_{ij}^2 ,  \hspace{1.5cm} |\theta_i - \theta_j| \le \overline\theta_{ij} \hspace{0.7cm} (i,j) \in \mathcal{L}.
\end{equation}


Here, the objective function \eqref{eq:obj R_OPF} minimizes the total real power generation cost subject to the following constraints: real and reactive power flow balance at bus $i$ \eqref{eq:real balance}--\eqref{eq:reactive balance}, real and reactive power flow  from $i$ to $j$ \eqref{eq:real flow}--\eqref{eq:reactive flow},  shunt susceptance selection for bus $i$ and tap ratio selection for line $(i,j)$ \eqref{eq:shunt and tap selection}, voltage magnitude bounds at bus $i$ \eqref{eq:voltage}, binary restrictions \eqref{eq:binary alpha and beta}--\eqref{eq:binary alpha and beta1}, reactive and active power output of generator $i$ \eqref{eq:real and reactive gen}, apparent flow and phase angle limit for each line $(i,j)$ \eqref{eq:flow and angle limit}.

\subsection{An Alternative MINLP}
In this section, we propose an alternative MINLP formulation of the ROPF problem motivated by \cite{kocuk2016strong}. Let us define a set of new decision variables $c_{ii}$, $c_{ij}$ and $s_{ij}$, respectively representing the quantities $ |V_i|^2$, $ |V_i||V_j| \cos(\theta_i - \theta_j)$ and $s_{ij} := -|V_i||V_j| \sin(\theta_i - \theta_j)$ for $i \in \mathcal{B}$ and $(i,j) \in \mathcal{L}$.
%
%
We denote the lower (upper) bounds of variables $c_{ii}, c_{ij}, s_{ij}$ as $\underline c_{ii}, \underline c_{ij},  \underline s_{ij}$ 
($\overline c_{ii}, \overline c_{ij},  \overline s_{ij}$) and set them as follows:
\begin{equation*}\begin{split}
&   \underline c_{ii} := \underline V_i^2, \  \overline c_{ii} := \overline 
V_i^2 \hspace{4.7cm} \ i \in \mathcal{B} \\
&  \underline c_{ij}:=\overline V_i \overline V_j \cos(\overline\theta_{ij}), \ \overline c_{ij}:=\overline V_i \overline V_j \hspace{2.2cm} (i,j) \in \mathcal{L}\\
&  \underline s_{ij}:=-\overline V_i \overline V_j \sin(\overline\theta_{ij}), \ \overline s_{ij}:=\overline V_i \overline V_j \sin(\overline\theta_{ij}) \hspace{.7cm} (i,j) \in \mathcal{L}.
\end{split}\end{equation*}

We will now discuss the constraints in the alternative formulation and their relations with the MINLP in Section~\ref{sec:MINLP}.
The updated  version of the real power flow balance 
constraint~\eqref{eq:real balance} is given as:
\begin{equation}\label{eq:real balance minlp}
p_i^g - p_i^d =  \ g_{ii} c_{ii} +  \sum_{j\in\delta(i)} p_{ij} \hspace{3.2cm} i \in \mathcal{B}.
\end{equation}


Since the variable $b_{ii}$ can be eliminated from the formulation by substituting $\sum_{k\in \mathcal{S}_i}  b_{ii}^k \alpha_{i}^k $, the reactive power flow  equation~\eqref{eq:reactive balance} is first rewritten as follows:
\begin{equation}\label{eq:reactive balance minlp_0}
q_i^g - q_i^d = - \left(\sum_{k\in \mathcal{S}_i}  b_{ii}^k \alpha_{i}^k\right) |V_i|^2 + \sum_{j\in\delta(i)} q_{ij}  \hspace{1.0cm} i \in \mathcal{B}.
\end{equation}
Then, we define a new variable $\Gamma_i^k :=  c_{ii} \alpha_i^k$ to linearize \eqref{eq:reactive balance minlp_0} and include additional constraints as follows:
\begin{equation}\label{eq:reactive balance minlp_}\begin{split}
q_i^g - q_i^d=   - \sum_{k\in\mathcal{S}_{i}} {b_{ii}^k} {\Gamma_i^k}  + \sum_{j\in\delta(i)} q_{ij}  \hspace{1.9cm} i \in \mathcal{B}  \\
 \underline c_{ii} \alpha_i^k  \le \Gamma_i^k \le   \overline c_{ii} \alpha_i^k , \quad    c_{ii}=   \sum_{k\in\mathcal{S}_{i}} {\Gamma_i^k} \hspace{1.6cm} i \in \mathcal{B}.
  \end{split}\end{equation}

%
%

We now update power flow  constraints using a similar procedure. In particular, we substitute $1 / \tau_{ij}$ with $\sum_{l\in \cT_{ij}}  \beta_{ij}^l / \tau_{ij}^l$ into constraints~\eqref{eq:real flow} and \eqref{eq:reactive flow}. After defining the new variables  $\bar \Phi_{ij}^l :=  c_{ii}  \beta_{ij}^l$, $\Phi_{ij}^l :=  c_{ij}  \beta_{ij}^l$ and  $\Psi_{ij}^l :=  s_{ij}  \beta_{ij}^l$,  we rewrite the real and reactive power flow constraints  \eqref{eq:real flow}--\eqref{eq:reactive flow} together with  other equations necessary for the linearization as follows:
%
%
%
\begin{equation}\label{eq:real flow minlp_}
\begin{split}
 &  p_{ij} = \sum_{l\in\mathcal{T}_{ij}} G_{ij} \left( \frac{c_{ii}  \beta_{ij}^l}{(\tau_{ij}^l )^2}  + \frac{\Phi_{ij}^l }{\tau_{ij}^l }\right) - B_{ij}\frac{\Psi_{ij}^l }{\tau_{ij}^l} \hspace{0.25cm} (i,j) \in \mathcal{L} 
\\
& q_{ij} = \sum_{l\in\mathcal{T}_{ij}} -B_{ij} \left( \frac{\bar \Phi_{ij}^l}{(\tau_{ij}^l )^2}  + \frac{\Phi_{ij}^l}{\tau_{ij}^l}\right)  - G_{ij} \frac{  \Psi_{ij}^l }{\tau_{ij}^l } \hspace{0.02cm} (i,j) \in \mathcal{L} 
\\
&   \underline c_{ii} \beta_{ij}^l  \le  \bar \Phi_{ij}^l  \le   \overline c_{ii} 
   \beta_{ij}^l , \ \ c_{ii}=  \sum_{l\in\mathcal{T}_{i,j}} {\bar \Phi_{ij}^l } \hspace{0.9cm}(i,j) \in \mathcal{L}
\\
&   \underline c_{ij} \beta_{ij}^l   \le   \Phi_{ij}^l  \le   \overline c_{ij}\beta_{ij}^l ,  \ c_{ij}=  \sum_{l\in\mathcal{T}_{i,j}} {  \Phi_{ij}^l } \hspace{0.9cm}(i,j) \in \mathcal{L}
\\
&   \underline s_{ij}\beta_{ij}^l   \le    \Psi_{ij}^l  \le   \overline s_{ij} \beta_{ij}^l ,\ s_{ij}= \sum_{l\in\mathcal{T}_{i,j}} {  \Psi_{ij}^l } \hspace{0.75cm} (i,j) \in \mathcal{L}.
   \end{split}
\end{equation}

We also update the constraint on voltage  magnitude 
bounds~\eqref{eq:voltage} as follows:
\begin{equation} \label{eq:voltage minlp}
   \underline V_i^2 \le c_{ii} \le \overline  V_i^2  \hspace{4.8cm} i \in \mathcal{B}.
\end{equation}

Finally, we define the following consistency constraints for each line $(i,j)$:
\begin{equation}\label{eq:consistency}
c_{ij}^2+s_{ij}^2  = c_{ii} c_{jj} \hspace{3.8cm} (i,j) \in \mathcal{L} 
\end{equation}
\begin{equation}\label{eq:consistency2}
(\Phi_{ij}^l)^2+(\Psi_{ij}^l)^2  = \bar \Phi_{ij}^l c_{jj} \hspace{2.7cm} (i,j) \in \mathcal{L}
\end{equation}
\begin{equation}\label{eq:Arctanconsistency}
\theta_j- \theta_i = \arctan (s_{ij}/c_{ij})  \hspace{2.6cm} (i,j) \in \mathcal{L}.
\end{equation}
Equation \eqref{eq:consistency} preserves the trigonometric relation between the variables $c_{ii}, c_{ij}$ and $s_{ij}$. If we multiply  \eqref{eq:consistency} by $\beta_{ij}^l$, we can get a similar condition for the variables $\bar \Phi_{ij}^l, \Phi_{ij}^l$ and $\Psi_{ij}^l $.

The alternative formulation minimizes \eqref{eq:obj R_OPF} subject to  constraints \eqref{eq:binary alpha and beta}--\eqref{eq:real balance minlp} and \eqref{eq:reactive balance minlp_}--\eqref{eq:Arctanconsistency}.

\subsection{MISOCP Relaxation}
The feasible region of the alternative MINLP formulation is  non-convex due to  constraints \eqref{eq:consistency}--\eqref{eq:Arctanconsistency}.
Let us relax these constraints as follows:
\begin{equation} \label{eq:consistency relaxed}\begin{split}
   c_{ij}^2+s_{ij}^2  \le c_{ii} c_{jj} \hspace{2.9cm} (i,j) \in \mathcal{L} \\
   (\Phi_{ij}^l)^2+(\Psi_{ij}^l)^2  \le \bar \Phi_{ij}^l c_{jj} \hspace{2.7cm} (i,j) \in \mathcal{L}.
\end{split}\end{equation}
Then, an MISOCP relaxation is obtained as \eqref{eq:obj R_OPF}, \eqref{eq:binary alpha and beta}--\eqref{eq:real balance minlp}, \eqref{eq:reactive balance minlp_}--\eqref{eq:voltage minlp} and  \eqref{eq:consistency relaxed}. 

To strengthen the MISOCP relaxation, we also consider an outer-approximation of the arctangent constraint \eqref{eq:Arctanconsistency}. This is achieved by the inclusion of  four hyperplanes as described in~\cite{kocuk2018matrix}. We will use the abbreviation \texttt{MISOCPA} to refer to this stronger relaxation. Additionally, we generate cutting planes for each cycle in the cycle basis using a method called \textit{SDP Separation}, more details  can be found in \cite{kocuk2016strong}. We denote this further improved relaxation as \texttt{MISOCPA+}.



%

\begin{table*}[t]
\centering
\caption{Computational results  (time is measured in seconds).}\label{tab:Cost Minimization}
\begin{tabular}{|c|rrrr|rrrr|}
\hline
& \multicolumn{4}{|c|}{\bf \texttt{TCR2}} & \multicolumn{4}{|c|}{\bf \texttt{MISOCPA+}}  \\
\hline
\bf Case & \bf LB & \bf Time & \bf UB & \bf \%Gap & \bf LB & \bf Time & \bf UB &  \bf \%Gap \\ \hline 
3lmbd & 5769.87 & 0.65 & 5812.64 & 0.74 & 5783.94 & 0.53 & 5812.64 & 0.49 \\ 
5pjm & 15313.38 & 0.72 & 17551.89 & 12.75 & 16395.73 & 0.22 & 17551.89 & 6.59 \\ 
30ieee & 205.19 & 1.13 & 205.64 & 0.22 & 205.15 & 1.38 & 205.25 & 0.05 \\ 
118ieee & 3695.39 & 4.66 & 3720.08 & 0.66 & 3684.68 & 12.89 & 3714.91 & 0.81 \\ 
\bf Average&  & \bf 1.79 &  & \bf 3.59 &  & \bf 3.75 &  & \bf 1.99 \\ \hline
3lmbd\_api & 363.00 & 0.66 & 367.74 & 1.29 & 362.92 & 0.56 & 367.74 & 1.31 \\ 
6ww\_api & 273.76 & 0.53 & 273.76 & 0.00 & 273.66 & 0.38 & 273.76 & 0.04 \\ 
14ieee\_api & 319.12 & 0.93 & 323.29 & 1.29 & 318.65 & 0.76 & 321.09 & 0.76 \\ 
30as\_api & 559.96 & 2.38 & 571.13 & 1.96 & 556.71 & 0.92 & 571.13 & 2.52 \\ 
30fsr\_api & 213.93 & 2.16 & 372.14 & 42.51 & 227.57 & 0.95 & 372.11 & 38.84 \\ 
39epri\_api & 7333.40 & 2.59 & 7466.25 & 1.78 & 7259.19 & 11.63 & 7480.45 & 2.96 \\
118ieee\_api & 5932.26 & 4.47 & 10258.47 & 42.17 & 5910.20 & 14.23 & 10158.62 & 41.82 \\ 
\bf Average&  & \bf 1.96 &  & \bf 13.00 &  & \bf \bf 4.20 &  & \bf 12.61 \\ \hline
3lmbd\_sad & 5831.07 & 0.57 & 5992.72 & 2.70 & 5867.46 & 0.53 & 5992.72 & 2.09 \\ 
4gs\_sad & 321.55 & 0.58 & 324.02 & 0.76 & 323.65 & 0.16 & 324.02 & 0.12 \\ 
5pjm\_sad & 25560.36 & 0.62 & 26423.33 & 3.27 & 26419.23 & 0.21 & 26423.32 & 0.02 \\ 
9wscc\_sad & 5521.49 & 0.54 & 5590.09 & 1.23 & 5589.54 & 0.20 & 5590.09 & 0.01 \\ 
29edin\_sad & 31173.80 & 3.19 & 46933.31 & 33.58 & 36270.50 & 2.77 & 45886.11 & 20.96 \\ 
30as\_sad & 903.09 & 2.32 & 914.44 & 1.24 & 906.96 & 1.11 & 914.44 & 0.82 \\ 
30ieee\_sad & 205.30 & 0.96 & 205.79 & 0.24 & 205.27 & 1.51 & 205.37 & 0.05 \\ 
118ieee\_sad & 3869.62 & 4.66 & 4323.91 & 10.51 & 4003.35 & 9.30 & 4258.72 & 6.00 \\ 
\bf Average&  & \bf 1.68 &  & \bf 6.69 &  & \bf 1.97 &  & \bf 3.76 \\ \hline
\bf Overall Average&  & \bf 1.81 &  & \bf 8.36 &  & \bf 3.17 &  & \bf 6.64 \\
 \hline
    \end{tabular}
\end{table*}

\noindent
\section{Computational Experiments}

\subsection{Algorithm}


We first solve the continuous relaxation of the \texttt{MISOCPA} formulation by relaxing the integrality of $\alpha_{i}^k$ and $\beta_{ij}^l$ variables. Then, for each cycle in the cycle basis, we use the SDP separation method to generate cutting planes to separate this continuous relaxation solution from the feasible region of the SDP relaxation of the cycle. The separation process is parallelized over cycles. We repeat this procedure five times consecutively. Then, we solve the final \texttt{MISOCPA+} relaxation to obtain a lower (LB) bound, and then fix the binary variables in the MINLP formulation to obtain an upper bound (UB) from the remaining NLP using a local solver. The optimality gap is computed as $\%\text{Gap}=100\times(1-\text{LB}/\text{UB})$.  

\subsection{Results}
We compare the percentage optimality gap and the computational time of the \texttt{MISOCPA+} approach  with the  publicly available implementation of TCR relaxation of Type 2 (\texttt{TCR2}) from~\cite{bingane2019tight}. All computational experiments have been carried out on a 64-bit desktop with Intel Core i7 CPU with 3.20GHz processor and 64 GB RAM. Our code is written in Python language using Spyder environment. The solvers Gurobi, IPOPT and MOSEK are used to solve the \texttt{MISOCPA+} relaxation, NLP and separation problems, respectively. 

For the computational experiments, we use the OPF instances from the NESTA library; typical operating conditions, congested operating conditions (API) and small angle difference conditions (SAD). We only consider difficult instances in which the SOCP optimality gap is more than 1\%~\cite{kocuk2016strong}.

The sets of the discrete values are determined as $b_{ii}^k \in \{ 0,1\}$ for $i\in\cS$ and $\tau_{ij}^l \in \{ 0.9,0.95,1,1.05,1.1\} $ for $(i,j)\in\cT$, which  represent the on/off status of the shunt susceptance and values of the tap ratio, respectively.

The results of our computational experiments are reported in Table~\ref{tab:Cost Minimization}.
We observe that  \texttt{MISOCPA+} has smaller optimality gap in fourteen out of nineteen instances, and better or same upper bound in eighteen of these instances (except 39epri\_api). If we compare the averages of optimality gap, \texttt{MISOCPA+} outperforms \texttt{TCR2} in all types of NESTA instances. \texttt{MISOCPA+}  has the best performance on SAD instances and  dominates \texttt{TCR2} in all of them. Overall, we  note that \texttt{MISOCPA+} relaxation has more accurate solutions with $ 6.64\%$  optimality gap, on average, than \texttt{TCR} with $ 8.36\%$. In terms of computational time, \texttt{MISOCPA+} is slower with 3.17 seconds,  on average, than \texttt{TCR2} with 1.81.

\section{Conclusion}
In this letter, we propose an MISOCP-based approach, namely \texttt{MISOCPA+}, to approximate globally optimal solutions of the ROPF problem. The accuracy and efficiency of this approach are compared with \texttt{TCR2} using difficult  OPF instances from the NESTA library. The computational results indicate that \texttt{MISCOPA+} is quite promising to solve any type of instances accurately, especially the ones with small angle conditions.

\bibliographystyle{plain}
\bibliography{references}

\end{document}